\documentclass{amsart}
\newtheorem{theorem}{Theorem}[section]
\newtheorem{lemma}[theorem]{Lemma}

\newtheorem{notation}[theorem]{Notation}
\newtheorem{corollary}[theorem]{Corollary}
\newtheorem{question}[theorem]{Question}

\theoremstyle{definition}

\theoremstyle{remark}
\newtheorem{remark}[theorem]{Remark}

\def \smash {\wedge}

\def \C {\mathbb{C}}

\def \bbH {\mathbb H}

\def \R {\mathbb{R}}
\def \T {\mathbb{T}}
\def \Z {\mathbb{Z}}

\def \CU {{\mathcal U}}

\begin{document}

\title[Almost complex $4$-manifolds]{Almost complex $4$-manifolds
  with vanishing first Chern class}
\author[Stefan Bauer]{Stefan Bauer} 
\address{ Fakult\"at f\"ur Mathematik,
Universit\"at Bielefeld, PF 100131, D-33501 Bielefeld}
%\curraddr{Institute for Advanced Study, 
%Einstein Drive, Princeton, NJ 08540}
\email{bauer@mathematik.uni-bielefeld.de}
%\thanks{This material is based upon work supported by the National
%Science Foundation under agreement No. DMS-0111298. Any opinions, findings and
%conclusions or recommendations expressed in this material are those of the author
%and do not necessarily reflect the views of the National Science Foundation.}
%\date{\today}
\begin{abstract} An odd Seiberg-Witten invariant imposes
        bounds on the signature of a closed, almost complex
        $4$-manifold with vanishing first Chern class. This
        applies in particular to symplectic $4$-manifolds of
        Kodaira dimension zero. 
\end{abstract}
\maketitle

%%%%%%%%%%%%%%%%%%%%%%%%%%%%%%%%%%%%%%%%%%%%%%%%%%%%%%%%%%%%%%%%%%%%%%%%%%%
%%%%%%%%%%%%%%%%%%%%%%%% Kapitel 0 %%%%%%%%%%%%%%%%%%%%%%%%%%%%%%%%%%%%%%%%
%%%%%%%%%%%%%%%%%%%%%%%%%%%%%%%%%%%%%%%%%%%%%%%%%%%%%%%%%%%%%%%%%%%%%%%%%%%

\section{Introduction\label{Intro}}

Vanishing of the first Chern class imposes severe
restrictions on a compact complex surface: It has to be
minimal and of Kodaira dimension at most zero. The
list of examples \cite{BPV}, p. 188, is rather short and known
to be complete   {(\cite{LYZ}, \cite{Teleman})}. It comprises in particular
$K3$-surfaces and tori, but also other examples found by and
named after Bombieri, Inoue, Hopf and Kodaira. Amongst these
surfaces, only the $K3$-surfaces exhibit nonvanishing
signature. 

Including closed symplectic $4$-manifolds into the
consideration, a few more examples of such with vanishing first Chern
class become available \cite{Thurston}, \cite{FGG},
\cite{Geiges}. However, $K3$-surfaces remain the
only known examples with nonvanishing signature.
The main result of this paper relates this more or less
empirical fact to Seiberg-Witten theory.

\begin{theorem}Let $X$ be a closed, almost complex
        $4$-manifold with vanishing first Chern class. If 
        the dimension $b_2^+(X)$ of a maximal positive
        definite linear subspace in the second cohomology of
        X satisfies $b_2^+(X)\geq 4$, then the Seiberg-Witten
        invariant of $X$ is an even number. 
\end{theorem} 

According to a theorem of Taubes \cite{Taubes94}, the absolute
value of the
Seiberg-Witten invariant of a symplectic $4$-manifold is
$1$, as soon as $b^+_2(X)\ge 2$. 
So this theorem applies, in particular, to compact
symplectic $4$-manifolds.

\begin{corollary}\label{symplectic} A closed,
        symplectic $4$-manifold $X$
        with torsion first Chern class satisfies the inequality
         $$b_2^+(X)\leq 3$$.          
\end{corollary}

Indeed, if the first Chern class of $X$ is torsion, then
there is a finite covering $\tilde X$ with vanishing first
Chern class, which of course is symplectic. 
The induced map $H^2(X;\R)\to
H^2({\tilde{X}};\R)$ is injective.

\begin{remark} Let $X$ be an almost complex manifold with
        vanishing first Chern class and $b_2^+(X)\leq 3$. 
        Then the signature is either zero or
        $-16$. If the signature is $-16$, then $b^+_2(X)=3$
        and the first
        Betti number vanishes. Otherwise the
        Betti numbers are related through $b_1=1+b_2^+$.
\end{remark}

Indeed, according to Rochlin's theorem the signature is
divisible by $16$. The assumption thus implies an upper bound $0$ for the
signature of $X$. In the equality $c_1^2-2c_2=p_1$ of
characteristic classes, the second Chern class equals the
Euler class of an almost complex manifold 
and the first Pontrjagin class
describes three times the signature. So we obtain an equality
$$0=\langle
2c_2+p_1,[X]\rangle=2(2-2b_1+b^+_2+b_2^-)+3(b_2^+-b_2^-)
=4(1-b_1+b_2^+)+{\rm{sign}}(X),
$$
from which we conclude $\rm{sign}(X)\ge -16$ and the claimed
values of the Betti numbers.

\begin{corollary}\label{negative} Let $X$ be a closed symplectic $4$-manifold
        with vanishing first Chern class and 
        $\rm {sign}(X)=-16$. Then the fundamental group of $X$ has
        no proper subgroup of finite index. 
\end{corollary}

Indeed, the covering manifold $\tilde{X}$ associated to a
subgroup of finite index $n$ would be compact symplectic with signature
${\rm{sign}}(\tilde{X})=n\cdot {\rm{sign}}(X)$ and with
$c_1({\tilde{X}})=0$.

Fundamental groups of complex surfaces are rather
restricted. In contrast, every finitely
presented group can be realized as the fundamental group of a symplectic
$4$-manifold \cite{Gompf}. But, doesn't any finitely
generated group contain subgroups of finite index? The amazing
answer is: No. In 1965, Richard Thompson 
constructed infinite simple groups which are
finitely presented, compare \cite{CFP}. 
If such a Thompson group $T$ did
admit a proper subgroup $S$ of
finite index, then the kernel of the action of $T$ on $T/S$
would be of finite index, contradicting simplicity. 

\begin{question} Is there a symplectic homology-$K3$-surface
with vanishing first Chern class and nontrivial fundamental
group? More specifically, can a Thompson group  be the
fundamental group of such a manifold?
\end{question}

The fundamental group of a symplectic manifold with
vanishing first Chern class and vanishing signature has a
property corresponding to \ref{negative}: Any subgroup of finite index has
rank at most $4$. Of course, this narrows the range of possible
fundamental groups of such manifolds. But still there is a considerable
gap if one compares with the
groups known to be realizable by symplectic manifolds of
Kodaira dimension zero.

Partial results with regard to \ref{symplectic} were obtained by
Morgan-Szabo \cite{MorganSzabo} under the assumption
$b_1(X)=0$ 
and by Tian-Jun Li \cite{Li} under the assumption
$b_1(X)\leq 4$. The main theorem above proves the 
"Betti Number Conjecture" in \cite{Li}.

The proof of the main result in the present
paper is modelled on the stable cohomotopy proof 
\cite{BauerSW}, thm 9.5, of Morgan-Szabo's result.
The concept can be explained in a few words: In its stable
homotopy interpretation  (\cite{Baueralt},
\cite{BauerFuruta}, \cite{BauerSW}), 
the Seiberg-Witten invariant is
the degree of a monopole map. Source and
target depend on index data of the given
$4$-manifold in a controlable way.
So it suffices to show that under the assumptions of
the theorem there are only maps of even degree between the
relevant spaces. This follows from 
equivariant obstruction theory using the fact that the vanishing
of the first Chern class leads to additional symmetry of the
monopole map.

\bigskip{\it Acknowledgement:} I am grateful to Katrin Tent
for helpful remarks on group theory.

\section{Proof of the main theorem}

Due to the vanishing of the integral first Chern class, the
almost complex $4$-manifold $X$ may be equipped with a spin
structure. Fixing such a spin structure and, furthermore, a
Riemannian metric on $X$ gives rise to a monopole map and a
refined Seiberg-Witten invariant. A key observation, well
known in gauge theory, is that for spin-manifolds the
monopole map is actually $Pin(2)$-equivariant, where
$Pin(2)\subset Sp(1) \subset {\bbH}^*$ is the normalizer
of the maximal torus ${\T}= Sp(1)\cap \C^*$ in the group 
$Sp(1)$ of quaternions of unit length.
As explained in \cite{BauerSW} ch. 9, the refined
invariant $[\mu]$ is a morphism in the $Pin(2)$-Spanier-Whitehead
category indexed by a universe $\CU$ containing only the
quaternions $\bbH$ and the real $1$-dimensional nontrivial 
$Pin(2)$-representation $V$ as irreducible summands.
The monopole morphism
$$ 
[\mu]\in 
\{T({\mathrm{ind}} D), S^{H^+(X)}\}^{Pin(2)}_{\CU}.
$$
is a morphism from the Thom spectrum of
the virtual index bundle of a 
family of Dirac operators to a sphere spectrum. Source and
target will be detailed in an instant.

Twisting the Dirac operator associated to the fixed spin
structure with flat $\T$-connections
defines a family of Dirac operators, parametrized by the
torus
$$Pic^0(X)=H^1(X;\R)/H^1(X;\Z).$$
The untwisted Dirac operator itself is $Sp(1)$-equivariant. 
This symmetry reduces to a $Pin(2)$-symmetry over the given
parameter space, with $j\in Pin(2)$ acting on $Pic^0(X)$ via
multiplication by $-1$. So the virtual index bundle can be 
represented as a difference 
$${\mathrm{ind}} D = F_0 - F_1$$
of complex vector bundles with quaternionic structures
(i.e. equipped with complex anti-linear bundle maps 
$j$ with $j^2=-1$ over the given involution of the base)
and thus an element in the Grothendieck group $KQ(Pic^0(X))$ of
such bundles.
Without loss of generality one can assume $F_1$ to be
trivial $F_1={\underline{\bbH}^c}:= Pic^0(X) \times
\bbH^c$. The rank of $F_0$ as
a complex bundle is determined by index theory:
$${\mathrm rk}_{\C}F_0 = \frac{-\rm{sign}(X)}{8}+2c.$$

The target of the monopole morphism is a sphere spectrum
$S^{H^+(X)}$, the suspension spectrum of the one-point
completed $b_2^+(X)=b$-dimensional vector 
space $H^+(X)$ of self-dual harmonic 
$2$-forms on $X$. The $Pin(2)$-action factors through the
quotient $\Z/2$ with $j$ acting by multiplication with $-1$; 
so we may identify $H^+(X)\cong V^{b}$. 

In particular, the monopole morphism is represented by
a $Pin(2)$-equivariant based map 
$$\mu:T(F_0)\smash S^{V^d}\to 
S^{\bbH^c + V^{b+d}}.$$
We are going to classify the relevant $Pin(2)$-equivariant maps from
Thom spaces to spheres as above using equivariant obstruction theory
as in \cite{tD}, ch. II.3. For this we need a slightly more
general setup:

\begin{notation}

\begin{enumerate}
\item\label{Euler pair}
Let $F$ be a complex vector bundle with quaternionic
structure over $Pic^0(X)$ and let $U\hookrightarrow W$ be a
fixed linear inclusion of $Pin(2)$-representations $U\cong V^{d_0}$ and 
$W\cong \bbH^c + V^{d_1}$.  We will call the pair 
$(T{\tilde F},S^W)$
consisting of the Thom space of the bundle $\tilde F= F+{\underline
  U}$ and the $1$-point
compactification $S^W$ of $W$ an {\rm Euler pair} of index
$(\xi,n)$ with $\xi=F-{\underline{\bbH}^c}\in KQ(Pic^0(X))$
and $n=d_1-d_0\ge 0$. The virtual dimension of the Euler
pair is the difference $\dim(T\tilde F)-\dim(S^W)=
b_1(X)+2\rm{rk}_\C(\xi)-n$ of the dimensions of the two
spaces.
\item\label{quasipole}
A map $f:T\tilde F\to S^W$ 
for an Euler pair $(T\tilde F, S^W)$ will be called a 
{\rm quasipole map}, if its 
restriction to the $\T$-fixed point set is the 1-point
        completion of the projection to the fiber
\[
Pic^0(X)\times U \to U \hookrightarrow W^\T.
\]
A quasipole map need not be 
$Pin(2)$-equivariant. 
The set of homotopy classes of $G$-equivariant quasipole
maps for subgroups $G\leq Pin(2)$ will be
denoted by $[T\tilde F,S^W]^G_q$.
\end{enumerate}
\end{notation}

The following lemma lists a few easy observations:

\begin{lemma}
Let $(T\tilde F, S^W)$ be an Euler pair.
\begin{enumerate}
\item The space $T\tilde F$ can be equipped with the
        structure of a
        $Pin(2)$-equivariant CW-complex. The space $S^W$ is
        a sphere with a linear $Pin(2)$-action.
%The isotropy groups of the $Pin(2)$-action on $T\tilde F$ and $S^W$ are
%$Pin(2)$, $\T$ and the trivial group.
%\item There are finitely many $Pin(2)$-fixed points in $TF$
%      and in $S^W$, respectively.
\item The $\T$-fixed point set  $T{\tilde F}^\T\subset
        T\tilde F$  is the Thom space $T\underline{U}$
        of the trivial bundle $\underline{U}= Pic^0(X)\times U$. 
        The residual $\Z/2=Pin(2)/\T$-action is via
        multiplication by $-1$ on both base and fiber of $\underline{U}$.
\item If the virtual dimension of $(T\tilde F,S^W)$ 
       is $1$ and $\dim(W)=w$, then the cohomology group 
$H^w(T\tilde F/\T,T\underline{U};\Z)$ is isomorphic
        to $\Z$.
\item 
Let $X$ be an almost complex manifold with
        vanishing first Chern class. Then the monopole morphism is 
        represented by a $Pin(2)$-equivariant quasipole map
        on an Euler pair of index $({\mathrm{ind}} D, b_2^+(X))$
        and virtual dimension $1$.
\end{enumerate}
\end{lemma}

\begin{proof}
After introducing an equivariant metric on the bundle
$F+\underline{U}$, the unit disc bundle is a manifold
with a differentiable $Pin(2)$-action and thus can be given
an equivariant CW-structure such that the sphere bundle is a
sub-complex. Such a CW-structure induces one on the Thom
space.

If an Euler pair is of virtual dimension $1$, then 
after replacing $T\underline{U}$ by a tubular neighborhood,
$H^*(T{\tilde F}/\T,T\underline{U};\Z)$  is the cohomology of  
a connected and orientable
manifold of dimension $w$ relative to its boundary.

The monopole morphism for a $4$-manifold is linear
when restricted to
the $\T$-fixed point sets and satisfies the defining
condition of a quasipole map, cf. \cite{Baueralt}, \cite{BauerSW}. 
For an almost complex manifold, the 
virtual dimension of an Euler pair with index
$({\mathrm{ind}} D, b^+_2(X))$ is (compare \ref{symplectic})
$$b_1(X)-\frac{\rm{sign}(X)}{4} - b_2^+(X)=1.$$
\end{proof}

In \cite{BauerSW} ch.4, a degree homomorphism 
$$h:
[T{\tilde F},S^W]^\T_q\to \Z$$
was defined for an Euler pair of index $(\xi,n)$ with $n\ge
2$, the sign depending on a choice of orientations.
The Seiberg-Witten invariant of $X$ is the degree of the
monopole morphism. The degree in the case of virtual
dimension $1$ is defined as follows:
The condition $n\ge 2$ implies a natural one-to-one correspondence
$$
[T{\tilde F}, S^W]_q^\T \cong [T{\tilde F}/T{\underline{U}}, S^W]^\T
$$
of $\T$-homotopy classes of quasipole maps with
$\T$-homotopy classes of maps which are constant on the 
$\T$-fixed point set. An element $f$ of the latter set induces a
homomorphism in reduced $\T$-equivariant Borel cohomology
$$f^*:{\tilde{H}}^*_\T(S^W;\Z)\to  
{{H}}^*_\T(T{\tilde F},T\underline{U};\Z)=
H^*(T\tilde F/\T,T\underline{U};\Z).$$
The image $f^*([W])\in H^*(T\tilde F/\T,T\underline{U};\Z)\cong\Z$ 
of the generator $[W]\in {\tilde{H}}^w_\T(S^W;\Z)$ of ${\tilde{H}}^*_\T(S^W;\Z)$ 
as a free $H^*_\T(pt)$-module is the degree
of $f$.

\begin{lemma}\label{Zwei} Let $(T\tilde{F},S^W)$ be an
      Euler pair of index $(\xi, n)$ and virtual dimension $1$. 
\begin{enumerate} 
\item There exists a $Pin(2)$-equivariant quasipole map
        $f:T{\tilde{F}}\to S^W$.
\item If $n\ge 2$, then the
degree map $h$
induces a bijection $$[T{\tilde{F}},S^W]_q^\T\cong \Z.$$
\item Suppose $n\geq 2$ and both $f$ and $g$ are
$Pin(2)$-equivariant quasipole maps on the given Euler pair.
Then the degrees of $f$ and $g$ differ by an even number.
\end{enumerate}
\end{lemma}

\begin{proof}
We have to show that 
the given map on the $\T$-fixed points
extends to a $Pin(2)$-equivariant map over $T\tilde F$.  
The obstructions
to extending over the 
$l$-skeleton of $T\tilde F$ are elements of obstruction groups
$${\mathfrak{H}}_{Pin(2)}^k(T\tilde{F}, T\underline{U};\pi_{k-1}(S^W))$$
(as defined in \cite{tD}, II.3)
for $k\leq l$. As long as $l\leq w$ holds, these
obstruction groups are zero due to the vanishing of
the coefficient groups. But the $w$-skeleton of $T\tilde{F}$ already
is the whole of $T{\tilde{F}}$ (a free $Pin(2)$-equivariant $k$-cell
$Pin(2)\times D^k$ has topological dimension $k+1$). This
proves the first part of the lemma.

If $n\ge2$, then we may use the natural one-to-one correspondence
$$[T{\tilde F}, S^W]_q^\T \cong [T{\tilde F}/T{\underline{U}}, S^W]^\T
$$
to single out a nullhomotopic $\T$-equivariant 
quasipole map.
Associating to a $\T$-equivariant quasipole map $f$ the difference cocycle to
the nullhomotopic quasipole map defines (thm II.3.17 in \cite{tD}) 
a bijective map 
$$[T{\tilde F},S^W]_q^\T\to 
{\mathfrak H}^w_{\T}(T{\tilde F},T\underline{U};\pi_w(S^W))\cong 
H^w(T\tilde F/\T,T\underline{U};\Z).$$
Composing with an isomorphism of the latter group with
$\Z$ results in the degree map $h$. 
This proves the second part of the lemma.

Obstruction theory associates to the  $Pin(2)$-equivariant 
quasipole maps $f$ and $g$ a difference cocycle and thus an element 
in the group ${\mathfrak
  H}^w_{Pin(2)}(T\tilde F,T\underline{U};\pi_w(S^W))$.
The image of the restriction homomorphism
$$ {\mathfrak
  H}^w_{Pin(2)}(T\tilde F,T\underline{U};\pi_w(S^W))\to
 {\mathfrak H}^w_{\T}(T\tilde F,T\underline{U};\pi_w(S^W))\cong \Z$$
consists of the multiples of $2= |Pin(2)/\T|$ by
\cite{tD}, prop. 4.9.
\end{proof}

The proof of the theorem uses the $Pin(2)$-equivariant Hopf map
$\eta: S^\bbH\to S^{V^3}$: 
Suppose the action of $Pin(2)$ on $\bbH$ is by
left multiplication and
consider $V^3\subset \bbH$ as embedded as purely imaginary
quaternions. Then the Hopf map is given by $\eta(h)=
\overline{h}ih$.

\medskip
We are now ready to finish the proof of the main theorem:
\begin{proof} 
Let $X$ be an almost complex manifold with vanishing first
Chern class and $b=b^+_2(X)\ge 4$. We will show that
the degree of every the $Pin(2)$-equivariant quasipole map
and thus the degree of the monopole morphism is even.
Because of the last part of
\ref{Zwei},
it suffices to exhibit a $Pin(2)$-equivariant 
quasipole map of even degree on an Euler pair of
index $({\mathrm{ind}} D,b)$. For this we
choose an Euler pair $(T\tilde F, S^W)$ of index 
$({\mathrm{ind}} D-\underline{\bbH},b-4)$. According to the first
part of \ref{Zwei}, there exists a $Pin(2)$-equivariant
  quasipole map  $f:T{\tilde{F}}\to S^W$. 
The map $f\smash \eta$ then is a $Pin(2)$-equivariant
quasipole map $T({\tilde F} +\underline{\bbH})\to
  S^{W+V^3}$. Composed
with the inclusion $S^{W+V^3}\hookrightarrow S^{W+V^4}$ we
get a $Pin(2)$-equivariant quasipole map of degree zero on
the Euler pair $(T(\tilde F+\underline{\bbH}), S^{W+V^4})$
of index $({\mathrm{ind}} D,b)$.
\end{proof}

%%%%%%%%%%%%%%%%%%%%%%%%%%%%%%%%%%%%% LITERATUR %%%%%%%%%%%%%%%%%%%%%%%%%%%%%%%

\end{document}